\documentclass[12pt, reqno]{amsproc}

\usepackage{amsmath,amscd,amsfonts,amssymb}
\usepackage[utf8]{inputenc}

\begin{document}
\title[Jordan Derivations of Special Subrings of Matrix Rings ]{ Jordan
Derivations of Special Subrings of Matrix Rings }
\author{ Umut Say\i n }
\address{ Department of Mathematics, Düzce University , 81620 Konuralp,Düzce, Turkey}
\email{umutsayin@duzce.edu.tr}
\author{Feride Kuzucuo\u{g}lu }
\address{ Department of Mathematics, Hacettepe University , 06800
Beytepe,Ankara, Turkey }
\email{feridek@hacettepe.edu.tr}
\maketitle

\begin{center}
\textbf{Abstract} 
\end{center}

Let $K$ be a 2-torsion free ring with identity and $R_{n}(K,J)$ be the ring
of all $n\times n$ matrices over $K$ such that the entries on and above the
main diagonal are elements of an ideal $J$ of $K.$ We describe all Jordan
derivations of the matrix ring $R_{n}(K,J)$ in this paper. The main
result states that every Jordan derivation $\Delta $ of $R_{n}(K,J)$ is of
the form $\Delta =D+\Omega $ where $D$ is a derivation of $R_{n}(K,J)$ and $\Omega $ is an extremal Jordan derivation of $R_{n}(K,J).$

\section{Introduction }

Let $K$ be an associative ring with identity and let $M_{n}(K)$ be the ring
of all $n\times n$ matrices over $K.$ Jordan multiplication is defined by $x\circ y=xy+yx$ for any $x,y\in M_{n}(K).$ If an additive map $\Delta $
of $M_{n}(K)$ satisfies $\Delta (x\circ y)=\Delta (x)\circ y+x\circ \Delta
(y)=\Delta (x)y+y\Delta (x)+x\Delta (y)+\Delta (y)x$ for all $x,y\in
M_{n}(K)$ , then $\Delta $ is called a Jordan derivation of $M_{n}(K).$
Derivations of the ring $M_{n}(K)$ give trivial Jordan derivations.
However, there are proper Jordan derivations which are not ring derivations
for some matrix rings (see \cite{Ghosseiri}, \cite{KoleMalt}, \cite{Kuzucuoglu}).

In 1957, Herstein (\cite{Herstein}) proved that every Jordan derivation of a prime ring
of characteristic not 2 is a derivation. This result was extended for
semi-prime rings (see \cite{Bresar}) and for certain triangular matrix algebras and
rings with nonzero nilpotent ideals (see \cite{Benkovic}, \cite{OuWangYao}, \cite{Zhang}).

Let $NT_{n}(K)$ be the ring of all $n\times n$ matrices over $K$ which are
all zeros on and above the main diagonal. Derivations and Jordan derivations
of the nilpotent ring $NT_{n}(K)$ were described in \cite{LevcRadc} and \cite{Kuzucuoglu}.

Let $M_{n}(J)$ be the ring of all $n\times n$ matrices over an ideal $J$ of $
K$ and $R=R_{n}(K,J)=NT_{n}(K)+M_{n}(J).$ All derivations of the matrix ring 
$R_{n}(K,J)$ were studied in \cite{KuzuSayi}.

Let $e_{i,j}$ denote the matrix with $1$ in the position $(i,j)$ and $0$
in every other position. The ring $R_{n}(K,J)$ is generated by the sets $
Ke_{i+1,i}$ ($i=1,2,...,n-1$) and $Je_{1,n}.$

The set of ideals $I=\{I_{i,j}:1\leq i,j\leq n\}$ of the ring $K$ are
called carpet if $I_{ij}I_{jk}\subseteq I_{ik}$ for any $i,j,k$ (see \cite{KargMerz}).
By using carpet ideals, we can easily compute powers of the ring $R.$ For
any matrix $[x_{i,j}]$ of the ring $R=R_{n}(K,J),$ the (i,j) entry $x_{i,j}$
is an element of $I_{i,j}$ where\ $I_{i,j}=K$ for $i>j$ and $I_{i,j}=J$ \
for $i\leq j.$ 

Let $Ann_{K}J=\{x\in K|\hspace*{0.1cm}xJ=Jx=0\}.$ Then $
AnnR=(Ann_{K}J)e_{n,1}$ $($see \cite{KuzuLevc}$).$ In this paper we determine the
structure of Jordan derivations on $R=R_{n}(K,J).$

\section{Construction of Standard Jordan Derivations of $R$}

(A1) Choose arbitrary additive group homomorphisms $\alpha ,\beta ,\gamma
:J\rightarrow Ann_{K}J$ satisfying $\alpha (J^{2})=0,$ $\beta (J^{2})=0$ $
,\gamma (J^{2})=0.$ Consider the following map of the set of all elementary
matrices

$\left. 
\begin{array}{c}
ye_{1,n}\rightarrow \alpha (y)e_{n-1,1}+\beta (y)e_{n-1,2}+\gamma (y)e_{n,2}
\\ 
ye_{1,n-1}\rightarrow \alpha (y)e_{n,1}+\beta (y)e_{n,2}\text{ \ \ \ \ \ \ \
\ \ \ \ \ \ \ \ \ \ \ \ \ \ \ } \\ 
ye_{2,n-1}\rightarrow \beta (y)e_{n,1}\text{ \ \ \ \ \ \ \ \ \ \ \ \ \ \ \ \
\ \ \ \ \ \ \ \ \ \ \ \ \ \ \ \ \ \ \ \ \ } \\ 
ye_{2,n}\rightarrow \beta (y)e_{n-1,1}+\gamma (y)e_{n,1}\text{ ,\ \ }
y\epsilon J\text{ \ \ \ \ \ \ \ \ }
\end{array}
\right\} $ \ \ \ \ \ (1)

We assume that images of the remaining elementary matrices from $R$ are
zeros.

If map (1) determines a Jordan derivation of the ring $R,$ then the
relations $ye_{1,n}\circ xe_{n,n-1}=yxe_{1,n-1},$ $xe_{2,1}\circ
ye_{1,n-1}=xye_{2,n-1},$ $xe_{2,1}\circ ye_{1,n}=xye_{2,n}$ give

$\left. 
\begin{array}{c}
\alpha (yx)=x\alpha (y) \\ 
\beta (yx)=x\beta (y) \\ 
\beta (xy)=\beta (y)x \\ 
\gamma (xy)=\gamma (y)x
\end{array}
\right\} $ \ \ \ \ \ \ \ \ \ \ \ \ \ \ \ \ \ \ \ \ \ \ \ \ \ \ \ \ \ \ \ \ \
\ \ \ \ \ \ \ \ \ \ \ \ \ \ \ (2)

for all $y\epsilon J$ and $x\epsilon K.$ One can easily check that if $
\alpha ,\beta ,\gamma $ satisfy (2), then (1) determines a proper Jordan
derivation of the ring $R$ but not a derivation unless $\alpha ,\beta
,\gamma $ are all zero maps. This Jordan derivation is called an 'extremal
Jordan derivation'.

\textbf{Example:} Let $K_{1}$ be a commutative ring with identity and let $J_{1}$ be
an ideal of $K_{1}$ which is nilpotent of class two. Let $K$ be a direct
product ($K_{1},K_{1}$) of two copies of the ring $K_{1}$ and let $
J=(J_{1},J_{1})$ be an ideal of $K.$ Then the maps
\begin{eqnarray*}
\alpha  &:&J\rightarrow Ann_{K}J \\
&:&(a,b)\rightarrow (a,0)
\end{eqnarray*}
\begin{eqnarray*}
\beta  &:&J\rightarrow Ann_{K}J \\
&:&(a,b)\rightarrow (0,b)
\end{eqnarray*}
\begin{eqnarray*}
\gamma  &:&J\rightarrow Ann_{K}J \\
&:&(a,b)\rightarrow (a,b)
\end{eqnarray*}
satisfy all the conditions of (2). Moreover, $\alpha (J^{2})=\beta
(J^{2})=\gamma (J^{2})=0.$

\bigskip (A2) For arbitrary additive group homomorphisms $\alpha _{1},\alpha
_{2}:J\rightarrow Ann_{K}J$ satisfying $\alpha _{1}(J^{2})=0,$ $\alpha
_{2}(J^{2})=0,$ the following map of the set $R_{3}(K,J)$ determines a
Jordan derivation of $R_{3}(K,J)$ under the conditions $\alpha
_{1}(xy)=\alpha _{1}(y)x,$ $\alpha _{2}(yx)=x\alpha _{2}(y)$ for $x\epsilon K
$ and $y\epsilon J.$ We assume that images of the elementary matrices except 
$ye_{1,3},$ $ye_{1,2}$ and $ye_{2,3}$ are zeros. 
\begin{eqnarray*}
ye_{1,3} &\rightarrow &\alpha _{1}(y)e_{3,2}+\alpha _{2}(y)e_{2,1} \\
ye_{1,2} &\rightarrow &\alpha _{2}(y)e_{3,1} \\
ye_{2,3} &\rightarrow &\alpha _{1}(y)e_{3,1}
\end{eqnarray*}

\bigskip (A3) If $\delta _{i}:J\rightarrow J,$ $\beta _{i}:J\rightarrow K,$ $
\theta :J\rightarrow K$ and $\gamma :J\rightarrow K$ $\ (i=1,2,3)$ are
additive mappings satisfying $\delta _{2}(J^{2})=0,$ $\beta _{2}(J)\subseteq
Ann_{K}J,$ $\delta _{1}(yz)=z\beta _{1}(y)+y\delta _{1}(z)+\delta _{1}(z)y,$ 
$y\beta _{3}(z)=\delta _{1}(yz),$ $\delta _{1}(yz)=y\theta (z),$ $\delta
_{2}(y)z+z\delta _{2}(y)=0,$ $z\gamma (y)+\theta (z)y=0,$ $\delta
_{3}(yz)=\beta _{1}(y)z,$ $\delta _{3}(yz)=\delta _{3}(y)z+z\delta
_{3}(y)+\beta _{3}(z)y,$ $\gamma (y)z=\delta _{3}(yz),$ $z\delta
_{1}(y)+\delta _{3}(y)z=\beta _{1}(yz)+\beta _{3}(yz),$ $\gamma (y)x=\beta
_{1}(yx)+\beta _{2}(xy),$ $\beta _{1}(yz+zy)=\beta _{1}(y)z+\beta _{1}(z)y,$ 
$z\beta _{1}(y)+\beta _{3}(z)y=0,$ $x\theta (y)=\beta _{2}(yx)+\beta
_{3}(xy),$ $\theta (yz)=y\theta (z)+z\beta _{2}(y),$ $\theta (y)z+y\beta
_{1}(z)=0,$ $\theta (xy)=x\delta _{1}(y)+\delta _{2}(y)x,$ $\theta
(yz)=y\beta _{3}(z),$ $z\gamma (y)+\beta _{3}(z)y=0,$ $\gamma (yz)=\beta
_{1}(y)z,$ $\gamma (yz)=\gamma (y)z,$ $\delta _{3}(y)x+x\delta
_{2}(y)=\gamma (yx),$ $z\gamma (y)+\delta _{1}(z)y+y\delta _{2}(z)=0,$ $
\delta _{1}(y)z+z\delta _{3}(y)+\delta _{1}(z)y+y\delta _{3}(z)=0,$ $\delta
_{2}(y)z+z\delta _{3}(y)+\theta (z)y=0,$ then the following map is a Jordan
derivation of $R_{3}(K,J)$ where $x\epsilon K,$ $y,z\epsilon J.$ We assume
that images of the elementary matrices except $ye_{1,3},$ $ye_{1,2}$, $
ye_{2,3}$ and $ye_{i,i}$ $(i=1,2,3)$ are zeros.
\begin{eqnarray*}
ye_{1,3} &\rightarrow &\sum\limits_{i=1}^{3}\delta _{i}(y)e_{i,i} \\
ye_{i,i} &\rightarrow &\beta _{i}(y)e_{3,1} \\
ye_{2,3} &\rightarrow &\theta (y)e_{2,1} \\
ye_{1,2} &\rightarrow &\gamma (y)e_{3,2}
\end{eqnarray*}

We now describe Jordan derivations which are also derivations of the ring $R$
defined in \cite{KuzuSayi}. 

\bigskip (B1) Inner Derivations: Let $A\epsilon R$. Then the derivation of $R
$ given by $X\rightarrow AX-XA,$ $X=[x_{i,j}]\epsilon R$ is called 'inner
derivation' induced by the matrix $A$.

\bigskip (B2) Diagonal Derivations: Let $D=\sum\limits_{i=1}^{n}d_{i}e_{i,i}$
($d_{i}\epsilon K$). Then the derivation of $R$ given by $X\rightarrow DX-XD,
$ $X=[x_{i,j}]\epsilon R$ is called 'diagonal derivation' induced by the
matrix $D.$

\bigskip (B3) Annihilator Derivations: An arbitrary 'annihilator derivation'
of the ring $R$ is of the form
\begin{equation*}
\lbrack a_{i,j}]\rightarrow \left[ \varsigma
_{n}(a_{1,n})+\sum\limits_{i=1}^{n-1}\varsigma _{i}(a_{i+1,i})\right] e_{n,1}
\text{ \ \ \ \ \ \ \ \ \ \ \ (3)}
\end{equation*}
where $\varsigma _{n}:J\rightarrow Ann_{K}J,$ $\varsigma _{n}(J^{2})=0$ and $
\varsigma _{i}:K\rightarrow Ann_{K}(J),$ $\varsigma _{i}(J)=0$ are additive
maps for $i=1,2,...,n-1.$

\bigskip (B4) Ring Derivations: If $\pi $ is a derivation of the coefficient
ring $K$ and the restriction of $\pi $ over the ideal $J$ is also a
derivation of $J,$ then $\pi $\ induces a derivation of the ring $R$ by the
rule
\begin{equation*}
\bar{\pi}(A)=\sum\limits_{i,j=1}^{n}\pi (a_{i,j})e_{i,j}\text{ , }
A=[a_{i,j}]\epsilon R
\end{equation*}

\bigskip (B5) Almost Annihilator derivation: If additive maps $\alpha ,\beta 
$:$J\rightarrow J,$ $\gamma :J\rightarrow K$ satisfy the relations $\alpha
(xy)=x\alpha (y),$ $\ \beta (yx)=\beta (y)x,$ $\ \ \ \ \gamma (y)z=y\gamma
(z)=\gamma (yz)=0$ and $\alpha (y)z+y\beta (z)=0,$ then the following map of
the set $R$
\begin{eqnarray*}
ye_{1,n} &\rightarrow &\alpha (y)e_{1,1}+\beta (y)e_{n,n}+\gamma (y)e_{n,1}
\\
ye_{i,n} &\rightarrow &\alpha (y)e_{i,1}\text{ \ , \ \ }1<i\leq n \\
ye_{1,j} &\rightarrow &\beta (y)e_{n,j}\text{ \ , \ \ }1\leq j<n
\end{eqnarray*}
determines a derivation of the ring $R$ called 'almost annihilator'
derivation.

\section{Jordan Derivations of $R_{n}(K,J)$}

Throughout this section, $K$ will be a 2-torsion free ring with identity.

\bigskip \textbf{Theorem 3.1.} Let $K$ be a 2-torsion free ring with identity, $J$ be
an ideal of $K$ and $R=R_{n}(K,J).$ If $n\geq 4$ then every Jordan
derivation of $R$ is of the form $\Delta =D+\Omega $ where $D$ is a
derivation of $R_{n}(K,J)$ and $\Omega $ is an extremal Jordan derivation
of $R_{n}(K,J).$ Moreover, $D$ is the sum of certain diagonal, inner,
annihilator, ring and almost annihilator derivations.

\bigskip 

For an arbitrary Jordan derivation $\Delta $ of the ring $R,$ we write the
image of an element $x_{i,j}e_{i,j}$ $\left( 1\leq i,j\leq n,\text{ }
x_{i,j}\epsilon I_{i,j}\right) $ of the form $\Delta
(x_{i,j}e_{i,j})=\sum\limits_{s,t=1}^{n}\Delta _{s,t}^{i,j}(x_{i,j})e_{s,t}$
where $\Delta _{s,t}^{i,j}$ are additive mappings from $I_{i,j}$ to $I_{s,t}$
.

\bigskip 

\textbf{Lemma 3.2.} Let $\Delta $ be an arbitrary Jordan derivation of R for $n\geq 4.
$ Then for $1<i<n-1$ and $x\epsilon K,$ $y\epsilon J$
\begin{equation*}
\Delta (xe_{i+1,i})=\sum \Delta
_{i+1,t}^{i+1,i}(x)e_{i+1,t}+\sum\limits_{s\neq i+1}\Delta
_{s,i}^{i+1,i}(x)e_{s,i}+\Delta _{n,1}^{i+1,i}(x)e_{n,1}\text{ \ \ \ \ (4)}
\end{equation*}
and
\begin{eqnarray*}
\Delta (ye_{1,n}) &=&\sum \Delta _{1,t}^{1,n}(y)e_{1,t}+\sum\limits_{s\neq
1}\Delta _{s,n}^{1,n}(y)e_{s,n}+\Delta _{n-1,1}^{1,n}(y)e_{n-1,1} \\
&&+\Delta _{n-1,2}^{1,n}(y)e_{n-1,2}+\Delta _{n,1}^{1,n}(y)e_{n,1}+\Delta
_{n,2}^{1,n}(y)e_{n,2}\text{ \ \ \ \ \ \ \ \ \ (5)}
\end{eqnarray*}

\bigskip 

\textbf{Proof} Let us fix $i,j$ and choose $k,m$ such that $k>m.$ If $k\neq j$ and $
m\neq i$ then $x_{i,j}e_{i,j}\circ y_{k,m}e_{k,m}=0$. By differentiating $
x_{i,j}e_{i,j}\circ y_{k,m}e_{k,m}=0,$ we get
\begin{eqnarray*}
0 &=&\Delta (x_{i,j}e_{i,j})\circ y_{k,m}e_{k,m}+x_{i,j}e_{i,j}\circ \Delta
(y_{k,m}e_{k,m}) \\
&=&\Delta (x_{i,j}e_{i,j})y_{k,m}e_{k,m}+y_{k,m}e_{k,m}\Delta
(x_{i,j}e_{i,j})+x_{i,j}e_{i,j}\Delta (y_{k,m}e_{k,m}) \\
&&+\Delta (y_{k,m}e_{k,m})x_{i,j}e_{i,j} \\
&=&\sum\limits_{s}\Delta
_{s,k}^{i,j}(x_{i,j})y_{k,m}e_{s,m}+\sum\limits_{t}y_{k,m}\Delta
_{m,t}^{i,j}(x_{i,j})e_{k,t} \\
&&+\sum\limits_{t}x_{i,j}\Delta
_{j,t}^{k,m}(y_{k,m})e_{i,t}+\sum\limits_{s}\Delta
_{s,i}^{k,m}(y_{k,m})x_{i,j}e_{s,j}
\end{eqnarray*}
Putting $y_{k,m}=1,$ the matrix on the right has zeros except $i-th$ , $k-th$
\ rows and $j-th$ , $m-th$ columns. Hence $\Delta _{s,k}^{i,j}=0$ for $m\neq
j,$ $s\neq i,$ $s\neq k$ and $\Delta _{m,t}^{i,j}=0$ for $i\neq k,$ $t\neq m,
$ $t\neq j.$ On the other hand, for $k>s>m$ and $k\neq i,j$ , $s\neq i,j$ , $
m\neq i,j$ , the (k,m), (s,m) and (k,s) coefficients of the equations $
\Delta (x_{i,j}e_{i,j}\circ e_{k,m})=0,$ $\Delta (x_{i,j}e_{i,j}\circ
e_{s,m})=0$ and $\Delta (x_{i,j}e_{i,j}\circ e_{k,s})=0$ are $\Delta
_{k,k}^{i,j}(x_{i,j})+\Delta _{m,m}^{i,j}(x_{i,j})=0...(\ast ),$ $\Delta
_{s,s}^{i,j}(x_{i,j})+\Delta _{m,m}^{i,j}(x_{i,j})=0...(\ast \ast )$ and $
\Delta _{k,k}^{i,j}(x_{i,j})+\Delta _{s,s}^{i,j}(x_{i,j})=0...(\ast \ast
\ast ),$ respectively.\ Comparing $(\ast ),(\ast \ast )$ and $(\ast \ast
\ast ),$\ we get $\Delta _{k,k}^{i,j}(x_{i,j})=\Delta
_{s,s}^{i,j}(x_{i,j})=\Delta _{m,m}^{i,j}(x_{i,j}).$\ By using $(\ast \ast
\ast ),$ $2\Delta _{k,k}^{i,j}(x_{i,j})=2\Delta
_{s,s}^{i,j}(x_{i,j})=0=2\Delta _{m,m}^{i,j}(x_{i,j}).$ Since K is 2-torsion
free, $\Delta _{k,k}^{i,j}(x_{i,j})=\Delta _{s,s}^{i,j}(x_{i,j})=\Delta
_{m,m}^{i,j}(x_{i,j})=0.$\ Therefore, the image of $xe_{i+1,i}$ ($x\epsilon K
$) under $\Delta $ is the matrix with zeros outside $i+1-th$ row, $i-th$
column and $(n,1)$ position. For $1<i<n-1,$ $\Delta (xe_{i+1,i})$ has the
form (4). In particular, $\Delta (xe_{2,1})=\sum \Delta
_{2,t}^{2,1}(x)e_{2,t}+\sum\limits_{s\neq 2}\Delta
_{s,1}^{2,1}(x)e_{s,1}+\Delta _{n,2}^{2,1}(x)e_{n,2}+\Delta
_{n,3}^{2,1}(x)e_{n,3}$ for $i=1$ and $\Delta (xe_{n,n-1})=\sum \Delta
_{n,t}^{n,n-1}(x)e_{n,t}+\sum\limits_{s\neq n}\Delta
_{s,n-1}^{n,n-1}(x)e_{s,n-1}+\Delta _{n-1,1}^{n,n-1}(x)e_{n-1,1}+\Delta
_{n-2,1}^{n,n-1}(x)e_{n-2,1}$ for $i=n-1.$ By the same relation $\Delta
_{n-1,1}^{1,n}(y)\neq 0$ since $k\neq n,$ $k\neq 1,$ $m\neq n-1$ and $\Delta
_{n-1,2}^{1,n}(y)\neq 0$ as$\ m\neq n-1,$ $k\neq 2$ for $y\epsilon J$.
Similarly, $\Delta _{n,1}^{1,n}(y)\neq 0$ since $k\neq 1,$ $m\neq n$ and $
\Delta _{n,2}^{1,n}(y)\neq 0$ as $m\neq n,$ $k\neq 2$ for $y\epsilon J$.
Thus we get (5).

\bigskip 

\textbf{Lemma 3.3.} Let $\Delta :R\rightarrow R$ be a Jordan derivation. Then there
exists a diagonal derivation $D$ of $R$ such that $(i+1,i)-th$ \ coefficient
of the matrix $\left( \Delta -\mathfrak{D}\right) (e_{i+1,i})$ is equal to
zero for all $i.$

\bigskip
\textbf{Proof} Let $D=\sum\limits_{i=1}^{n}d_{i}e_{i,i}$ where $d_{1}=0,$ $
d_{i+1}=\sum\limits_{k=1}^{i}a_{k}$ and $a_{k}=\Delta _{k+1,k}^{k+1,k}(1).$
Then there exists a diagonal derivation $D:X\rightarrow DX-XD$ by the
diagonal matrix $D$ such that $D(e_{i+1,i})=De_{i+1,i}-e_{i+1,i}D=\Delta
_{i+1,i}^{i+1,i}(1)e_{i+1,i}$ . Since $(i+1,i)-th$ coefficient of the matrix 
$\Delta $ is equal to $\Delta _{i+1,i}^{i+1,i}(1),$ the proof is completed.

\bigskip 
\textbf{Lemma 3.4.} Let $\Delta :R\rightarrow R$ be a Jordan derivation
such that $\left( i+1,i\right) -th$ coefficient of matrices $\Delta
(e_{i+1,i})$ are zero for all $1\leq i<n.$ Then there exists an inner
derivation $I$ such that each matrix $\left( \Delta -\mathfrak{I}\right)
(e_{i+1,i})$ has zero $i-th$ column and $(i+1,1)$ entries.

\bigskip \textbf{Proof} Define a matrix $A=[A_{i,j}]$ such that $A_{u,u}=0$ ($1\leq
u\leq n),$ $A_{u,i+1}=\Delta _{u,i}^{i+1,i}(1)$ ($u\neq i+1,$ $i\neq n$) and 
$A_{j,1}=0$ ($1\leq j\leq n$). Consider the action of $I_{A}$ on the
matrices $e_{i+1,i}$ . Then $I_{A}(e_{i+1,i})=Ae_{i+1,i}-e_{i+1,i}A=\sum
\limits_{\substack{ k=1 \\ k\neq i+1}}^{n}\Delta
_{k,i}^{i+1,i}(1)e_{k,i}+\sum\limits_{\substack{ m=1 \\ m\neq i-1}}
^{n-1}-\Delta _{i,m}^{m+1,m}(1)e_{i+1,m+1}$ . Therefore, i-th column of each
matrix $(\Delta -I_{A})(e_{i+1,i})$ is equal to zero. Now define a matrix $B=[B_{i,j}]$ such that $B=-b_{3}e_{2,1}-b_{4}e_{3,1}-...-b_{n}e_{n-1,1}$
where $b_{i+1}=\Delta _{i+1,1}^{i+1,i}(1) \ (1<i<n)$. Denote by $I_{B}$
the inner derivation induced by the matrix $B.$ It can be easily seen that $I_{B}(e_{i+1,i})=Be_{i+1,i}-e_{i+1,i}B=\Delta _{i+1,1}^{i+1,i}(1)e_{i+1,1}$ ($i=2,...,n-1$). Hence i-th columns and (i+1,1) entries of the matrices ($\Delta -I$)($e_{i+1,i}$) are zeros with $I=I_{A}+I_{B}.$

\bigskip \textbf{Lemma 3.5.} Let $\Delta :R\rightarrow R$ be a Jordan derivation such
that $i-th$ columns and $(i+1,1)$ entries of the matrices $\Delta (e_{i+1,i})
$ are all zeros. Then the following equalities are obtained for $
x_{k,m}\epsilon I_{k,m},$ $x\epsilon K$ and $y\epsilon J$ $;$
\begin{eqnarray}
\Delta (e_{2,1}) &=&0  \\
\Delta (e_{n,n-1}) &=&0 \\
\Delta (e_{i+1,i}) &=&\Delta _{n,1}^{i+1,i}(1)e_{n,1}\text{ , }1<i<n-1 \\
\Delta (e_{i,j}) &=&0\text{ , }i-j>1 \\
\Delta (xe_{i,1}) &=&\Delta _{i,1}^{i,1}(x)e_{i,1}+\Delta
_{n,1}^{i,1}(x)e_{n,1}\text{ \ \ , \ \ }1<i<n \\
\Delta (x_{n,j}e_{n,j}) &=&\Delta _{n,1}^{n,j}(x_{n,j})e_{n,1}+\Delta
_{n,j}^{n,j}(x_{n,j})e_{n,j} \\
\Delta (ye_{1,i}) &=&\Delta _{1,1}^{1,i}(y)e_{1,1}+\Delta
_{1,2}^{1,i}(y)e_{1,2} \notag \\
&&+\Delta _{1,i}^{1,i}(y)e_{1,i}+\Delta _{n,1}^{1,i}(y)e_{n,1} \notag \\
&&+\Delta _{n,2}^{1,i}(y)e_{n,2}+\Delta _{n,i}^{1,i}(y)e_{n,i} \text{  \ \ , \ \   } i\neq 1,n \\
\Delta (ye_{i,n}) &=&\Delta _{i,1}^{i,n}(y)e_{i,1}+\Delta
_{n-1,1}^{i,n}(y)e_{n-1,1} \notag \\
&&+\Delta _{n,1}^{i,n}(y)e_{n,1}+\Delta _{i,n}^{i,n}(y)e_{i,n} \notag \\
&&+\Delta _{n-1,n}^{i,n}(y)e_{n-1,n}+\Delta _{n,n}^{i,n}(y)e_{n,n} \\
\Delta (ye_{1,n}) &=&\Delta _{1,1}^{1,n}(y)e_{1,1}+\Delta
_{1,2}^{1,n}(y)e_{1,2}+\Delta _{1,n}^{1,n}(y)e_{1,n} \notag \\
&&+\Delta _{n-1,1}^{1,n}(y)e_{n-1,1}+\Delta
_{n-1,2}^{1,n}(y)e_{n-1,2}+\Delta _{n-1,n}^{1,n}(y)e_{n-1,n} \notag \\
&&+\Delta _{n,1}^{1,n}(y)e_{n,1}+\Delta _{n,2}^{1,n}(y)e_{n,2}+\Delta
_{n,n}^{1,n}(y)e_{n,n} \\
\Delta (x_{i,j}e_{i,j}) &=&\Delta _{i,1}^{i,j}(x_{i,j})e_{i,1}+\Delta
_{i,j}^{i,j}(x_{i,j})e_{i,j} \notag \\
&&+\Delta _{n,1}^{i,j}(x_{i,j})e_{n,1}+\Delta _{n,j}^{i,j}(x_{i,j})e_{n,j}
\text{ \ \ , \ \ }1<i,j<n
\end{eqnarray}

\bigskip \textbf{Proof} The Equations (1)-(10) can easily be obtained by the
relations below when $K$ is 2-torsion free;
\begin{eqnarray*}
\Delta \left( e_{i+1,i}\circ e_{i+1,i}\right)  &=&0, \\
\Delta (e_{i+1,i}\circ e_{j+1,j}) &=&0\text{ }(i\neq j-1,i\neq j,i\neq j+1),
\\
\Delta (e_{i+1,i}\circ e_{i+2,i}) &=&0, \\
\Delta (e_{i+1,i}\circ e_{i+1,i}) &=&0,\text{ }\Delta (e_{n,n-1}\circ
e_{n,n-1})=0, \\
\Delta (e_{n,n-1}\circ e_{n-1,n-2}) &=&\Delta (e_{n,n-2}),\text{ }\Delta
(e_{n,n-1}\circ e_{n,n-2})=0, \\
\Delta (e_{2,1}\circ e_{2,1}) &=&0,\text{ }\Delta (e_{2,1}\circ e_{3,1})=0,
\text{ }\Delta (e_{3,2}\circ e_{2,1})=e_{3,1}, \\
\Delta (x_{i,1}e_{i,1}\circ e_{k,m}) &=&0\text{ }(1<i<n,m\neq i,k>m), \\
\Delta (x_{n,i}e_{n,i}\circ e_{k,m}) &=&0\text{ }(i\neq 1), \\
\Delta (y_{1,i}e_{1,i}\circ e_{k,m}) &=&0\text{ }(k\neq i,m\neq 1,i\neq n),
\\
\Delta (y_{i,n}e_{i,n}\circ e_{k,m}) &=&0\text{ }(i\neq 1,i\neq n,k\neq
n,m\neq i), \\
\Delta (y_{1,n}e_{1,n}\circ e_{k,m}) &=&0\text{ }(k\neq n,m\neq 1), \\
\Delta (x_{i,j}e_{i,j}\circ e_{k,m}) &=&0\ (1<i,j<n).
\end{eqnarray*}

\textbf{Lemma 3.6.} Let $\Delta $ be a Jordan derivation of R satisfying the
conditions (1)-(10) in Lemma 3.5. Then there exists an annihilator
derivation $\Upsilon $ such that $(\Delta -\Upsilon )$($e_{i+1,i}$) is
equal to zero for all $i$ .

\bigskip \textbf{Proof} Let $x\epsilon K,$ $y,z\epsilon J$ be arbitrary elements. For 
$i\neq 1,n$ $,$ $\Delta (e_{n,i}\circ ye_{1,n})=\Delta (ye_{1,i})$. This
implies that $\Delta _{1,1}^{1,i}=0$ . If we combine $\Delta (e_{n,1}\circ
ye_{1,n-1})=\Delta (ye_{n,n-1})$ and $e_{n,1}\circ \Delta
(ye_{1,n-1})=\Delta (ye_{n,n-1})$ , we obtain $\Delta
_{11}^{1,n-1}(y)=\Delta _{n,1}^{n,n-1}(y)=0$ . Let us say $\varsigma
_{i}=\Delta _{n,1}^{i+1,i}$ . Then $\varsigma _{n-1}(J)=0.$\ If $k>2$ , $\Delta (e_{2,1}\circ ye_{1,k})=\Delta (ye_{2,k})$ follows that $\Delta
_{n,k}^{2,k}=0$. Again for $k>2$ , $\Delta (ye_{2,k}\circ e_{k,1})=\Delta
(ye_{2,1})$ and $\Delta (ye_{2,k})\circ e_{k,1}=\Delta (ye_{2,1})$ implies that $\Delta _{n,k}^{2,k}(y)=\Delta _{n,1}^{2,1}(y)=0=\varsigma _{1}(y)$.
For $i\neq 1,n-1$ , $\Delta (e_{n,i}\circ ye_{i+1,n})=\Delta (ye_{i+1,i})$
and $e_{n,i}\circ \Delta (ye_{i+1,n})=\Delta (ye_{i+1,i})$ follows $\Delta
_{i,1}^{i+1,n}(y)=\Delta _{n,1}^{i+1,i}(y)=0=\varsigma _{i}(y)$. So $\varsigma _{i}(J)$ is zero ($i<n$). For $2\leq i\leq n-2$ , we get $\varsigma _{i}=\Delta _{n,1}^{i+1,i}:K\rightarrow Ann_{K}J$ \ by using the
relations $\Delta (xe_{i+1,i}\circ ye_{1,k})=0$ ($k\neq i+1,n$)\ and $\Delta
(xe_{i+1,i}\circ ye_{k,n})=0$ ($k\neq 1,i$)\ . By $\Delta (xe_{2,1}\circ
ye_{n-1,n})=0$ and $\Delta (xe_{2,1}\circ ye_{1,k})=\Delta (xye_{2,k})$ \ ($k\neq 2,n$), we have $y\Delta _{n,1}^{2,1}(x)=0$ , $\Delta
_{n,1}^{2,1}(x)y=0$ and $\varsigma _{1}=\Delta _{n,1}^{2,1}:K\rightarrow
Ann_{K}J$. Similarly, by $\Delta (xe_{n,n-1}\circ ye_{1,2})=0$ and $\Delta (xe_{n,n-1}\circ ye_{1,n})=\Delta (yxe_{1,n-1})$ we obtain $\Delta
_{n,1}^{n,n-1}(x)y=0$ , $y\Delta _{n,1}^{n,n-1}(x)=0$\ and $\varsigma
_{n-1}=\Delta _{n,1}^{n,n-1}:K\rightarrow Ann_{K}J$. Now consider the product $\Delta (ye_{1,k}\circ ze_{k,n})=\Delta (yze_{1,n}).$ Hence $\Delta _{n,1}^{1,n}(yz)=0.$ Say $\varsigma _{n}=\Delta _{n,1}^{1,n}$ ,
then $\varsigma _{n}(J^{2})=0$. By $\Delta (ye_{1,1}\circ ze_{1,n})=\Delta
(yze_{1,n})$ and $\Delta (ye_{1,n}\circ ze_{n,n})=\Delta (yze_{1,n})$ , we
have $0=\Delta _{n,1}^{1,n}(z)y=z\Delta _{n,1}^{1,n}(y)$\ and $\varsigma
_{n}:J\rightarrow Ann_{K}J$ . Thus $\Upsilon :[x_{i,j}]\rightarrow \left(
\varsigma _{n}(x_{1,n})+\sum\limits_{i=1}^{n-1}\varsigma
_{i}(x_{i+1,i})\right) e_{n,1}$ is an annihilator derivation and $(\Delta
-\Upsilon )(e_{i+1,i})=0$ \ for all $i.$ Say $\Theta =\Delta -\Upsilon .$
Hence (n,1)-coefficients of $\Theta (xe_{i+1,i})$ and $\Theta (ye_{1,n})$
are equal to zero.

\bigskip 

\textbf{Lemma 3.7.} Let $\Theta =\Delta -\Upsilon $ be a Jordan derivation of the
ring R as in Lemma 3.6. Then there exists a ring derivation $\bar{\theta}$
of $R$ such that $(i,j)$ coefficient of $\Theta (x_{i,j}e_{i,j})$ is equal
to zero$.$

\bigskip \textbf{Proof} Let $x,x_{1},x_{2}\epsilon K$ and $y\epsilon J$ be arbitrary
elements. By using the relation $x_{1}e_{i,j}\circ
x_{2}e_{j,k}=x_{1}x_{2}e_{i,k}$ for $i>j>k$ we have $\Theta
_{i,j}^{i,j}(x_{1})x_{2}+x_{1}\Theta _{i,j}^{i,j}(x_{2})=\Theta
_{i,j}^{i,j}(x_{1}x_{2})$. If we say $\Theta _{i,j}^{i,j}=\theta $ , $
\theta $ is a ring derivation of K. Similarly, $\Theta
_{i,j}^{i,j}(y)=\Theta _{i,i-1}^{i,i-1}(y)$ by $ye_{i,j}\circ
e_{j,i-1}=ye_{i,i-1}$ for $i\leq j$. This implies that $\theta $ is a
derivation of J as well. So $\bar{\theta}:[x_{i,j}]\rightarrow
\sum\limits_{i,j}\theta (x_{i,j})e_{i,j}$ is a ring derivation of R.

\bigskip Let us say $\Xi =\Theta -\bar{\theta}$. Thus (i,j)-coefficients of
the matrices $\Xi (x_{i,j}e_{i,j})$ are equal to zero.

\bigskip \textbf{Lemma 3.8.} Let $\Xi $ be a Jordan derivation of R as in Lemma 3.7.
Then $\Xi (xe_{i,j})=0$ for all $i>j$ , $\Xi (ye_{i,i})=\Xi
_{n,1}^{i,i}(y)e_{n,1}$ for all $i,$ $\Xi (ye_{1,j})=\Xi
_{n,j}^{1,j}(y)e_{n,j}$ for $1<j<n-1$ , $\Xi (ye_{i,n})=\Xi
_{i,1}^{i,n}(y)e_{i,1}$ for $i\neq 1,2$ where $x\epsilon K,$ $y\epsilon J$
are arbitrary elements.

\bigskip \textbf{Proof} Let $x,x_{1},x_{2}\epsilon K$ and $y\epsilon J$ be arbitrary
elements. For $1<i<n$, by combining $\Xi (x_{1}e_{i+1,i}\circ
x_{2}e_{i,i-1})=\Xi (x_{1}x_{2}e_{i+1,i-1})$ and $\Xi (x_{1}e_{i+1,i})\circ
x_{2}e_{i,i-1}+x_{1}e_{i+1,i}\Xi (x_{2}e_{i,i-1})=\Xi
(x_{1}x_{2}e_{i+1,i-1})\ $it can be easily seen that $\Xi
_{n,i-1}^{i+1,i-1}=0=\Xi _{i+1,1}^{i+1,i-1}.$ Then $\Xi
(x_{1}x_{2}e_{i+1,i-1})=0$ . By using $\Xi (x_{1}e_{i+1,i}\circ
x_{2}e_{i,i-1})=\Xi (x_{1}x_{2}e_{i+1,i-1})=0$ , we get $\Xi
_{n,i}^{i+1,i}=0$, $\Xi _{i,1}^{i,i-1}=0$. Thus $\Xi (xe_{i,j})=0$\ for $i>j$. To obtain $\Xi (ye_{i,i})=\Xi _{n,1}^{i,i}(y)e_{n,1}$ , we use the
following relations; $\Xi (xe_{i+2,i+1}\circ ye_{i+1,i+1})=\Xi
(xye_{i+2,i+1})=0$ ($1\leq i<n-1$), $\Xi (xe_{2,1}\circ ye_{1,1})=\Xi
(xye_{2,1})=0$ and $\Xi (xe_{i+1,i}\circ ye_{i+1,i+1})=\Xi (yxe_{i+1,i})=0$
($1\leq i<n$). By $\Xi (xe_{i,1}\circ ye_{1,j})=\Xi (xye_{i,j})=0$ for $i>j$ and $\Xi (xe_{n,j}\circ ye_{1,n})=\Xi (yxe_{1,j})$ for $1<j<n-1$ , we
have $\Xi (ye_{1,j})=\Xi _{n,j}^{1,j}(y)e_{n,j}$ for $1<j<n-1$. Finally
by $\Xi (xe_{i,1}\circ ye_{1,n})=\Xi (xye_{i,n})$ for $1<i<n$ we have $\Xi
(ye_{i,n})=\Xi _{i,1}^{i,n}(y)e_{i,1}$.

\bigskip 

\bigskip 

\textbf{Lemma 3.9.} Let $\Xi $ be a Jordan derivation of R as in Lemma 3.8. Then
there exists an almost annihilator derivation $\Gamma $ of R such that $\Xi
-\Gamma $ is an extremal Jordan derivation of R which is defined in Section
2.

\bigskip \textbf{Proof} Let $\bar{\alpha}=\Xi _{1,1}^{1,n},$ $\bar{\beta}=\Xi
_{n,n}^{1,n},$ $x\epsilon K$ and $y,z\epsilon J$. By using the
corresponding appropriate relations and considering K as 2-torsion free, the
additive map $\Gamma $ on $R$ defined by $\Gamma (ye_{1,n})=\bar{\alpha}
(y)e_{1,1}+\bar{\beta}(y)e_{n,n}$ , $\Gamma (ye_{i,n})=\bar{\alpha}(y)e_{i,1}
$ ($1<i\leq n$), $\Gamma (ye_{1,j})=\bar{\beta}(y)e_{n,j}$ $(1\leq j<n)$, $\Gamma (x_{i,j}e_{i,j})=0$ $(1<i$ and $j<n)$ is obviously an almost
annihilator derivation of R.\ 

Consider the relations $\Xi (xe_{n-1,j}\circ ye_{1,n-1})=\Xi (yxe_{1,j})$ \
for $1<j<n-1$ and $\Xi (xe_{i,1}\circ ye_{1,j})=\Xi (xye_{i,j})$, $\Xi
(ye_{1,j})\circ xe_{i,1}=\Xi (xye_{i,j})$ for $1<i<j<n$ where $(i,j)\neq
(2,n-1),$ then we obtain $\Xi (ye_{1,j})=0$ for $j\neq n,n-1$ and $\Xi
(ye_{i,j})=0$ \ for $1<i<j<n$ respectively. 

Finally, if we say $\Pi =\Xi -\Gamma $ then$\ \Pi (ye_{1,n})=\Pi
_{n-1,1}^{1,n}(y)e_{n-1,1}+\Pi _{n-1,2}^{1,n}(y)e_{n-1,2}+\Pi
_{n,2}^{1,n}(y)e_{n,2}$ , $\Pi (ye_{1,n-1})=\Pi
_{n,1}^{1,n-1}(y)e_{n,1}+\Pi _{n,2}^{1,n-1}(y)e_{n,2}$, $\Pi
(ye_{2,n-1})=\Pi _{n,1}^{2,n-1}(y)e_{n,1}$ , $\Pi (ye_{2,n})=\Pi
_{n-1,1}^{2,n}(y)e_{n-1,1}+\Pi _{n,1}^{2,n}(y)e_{n,1}$ and $\Pi
(x_{i,j}e_{i,j})=0$ for $(i,j)\neq (1,n),(2,n),(1,n-1),(2,n-1)$. By using
relations $ye_{1,n}\circ xe_{n,n-1}=yxe_{1,n-1}$ , $xe_{2,1}\circ
ye_{1,n-1}=xye_{2,n-1}$ , $xe_{2,1}\circ ye_{1,n}=xye_{2,n}$ , we get $\alpha =\Pi _{n-1,1}^{1,n}=\Pi _{n,1}^{1,n-1}$ , $\beta =\Pi
_{n-1,2}^{1,n}=\Pi _{n,2}^{1,n-1}=\Pi _{n,1}^{2,n-1}=\Pi _{n-1,1}^{2,n}$ , $\gamma =\Pi _{n,2}^{1,n}=\Pi _{n,1}^{2,n}$\ and $\alpha (yx)=x\alpha (y),$ $
\beta (yx)=x\beta (y),$ $\beta (xy)=\beta (y)x,$ $\gamma (xy)=\gamma (y)x$.
Then for $n>4,$ the relations $0=\Pi (ye_{1,3}\circ ze_{3,n})=\Pi
(yze_{1,n})$, $0=\Pi (yze_{1,n})=\Pi (ye_{1,n-1}\circ ze_{n-1,n})=\Pi
(ye_{1,1}\circ ze_{1,n})$, $0=\Pi (yze_{1,n})=\Pi (ye_{1,2}\circ
ze_{2,n})=\Pi (ye_{1,n}\circ ze_{n-1,n-1})$ and $0=\Pi (ye_{1,n}\circ
ze_{3,n})$ gives\ all conditions of extremal Jordan derivation defined in
A1. 

\bigskip If $n=4$ , we use the following relations to show all conditions
of extremal Jordan derivation of $R_{4}(K,J)$ are satisfied;

\bigskip By $ye_{1,3}\circ ze_{3,4}=yze_{1,4}$ , $ye_{1,2}\circ
ze_{2,4}=yze_{1,4}$ , $ye_{1,1}\circ ze_{1,4}=yze_{1,4}$ we get $\alpha
(J^{2})=0,$ $\beta (J^{2})=0,$ $\gamma (J^{2}),$ respectively and by $
ye_{1,1}\circ ze_{1,4}=ye_{1,2}\circ ze_{2,4}=ye_{1,3}\circ
ze_{3,4}=yze_{1,4}$ , $ye_{1,4}\circ ze_{3,4}=0$ . 

\bigskip Now \textbf{Theorem 3.1} follows easily by the Lemmas 3.2 - 3.9.

\bigskip 

\bigskip For \textbf{n=3}, after applying Lemmas 3.2-3.7, it is obtained that $\Xi $
is equal to the sum of the Jordan derivations A2 and A3 described in
Section 2.

\end{document}